\documentclass[12pt]{article}
\usepackage{amsfonts, amsmath,  graphicx,  graphics,  subfigure,    mathrsfs}
\usepackage[left=3 cm,top=3cm,right=3 cm, bottom = 4cm]{geometry}
 
\newcommand{\naturals}{\mathbb{N}}
\newcommand{\reals}{\mathbb{R}}
\newcommand{\pr}{\mathbb{P}}
\newcommand{\ex}{\mathbb{E}}
\newcommand{\pspace}{(T,\mathscr{T},\mu)}
\newtheorem{theorem}{Theorem}

\newenvironment{example}[1][Example]{\begin{trivlist}
\item[\hskip \labelsep {\bfseries #1}]}{\end{trivlist}}

\newenvironment{remark}[1][Remark]{\begin{trivlist}
\item[\hskip \labelsep {\bfseries #1}]}{\end{trivlist}}

\begin{document}

\title{The Uniform Law for  Sojourn Measures of\\  Random Fields }

\author{K.\ Borovkov\footnote{
Department of Mathematics and Statistics, University of Melbourne, Parkville 3010, Australia. E-mail: borovkov@unimelb.edu.au.}
\ and \ S.\ McKinlay
\footnote{
Department of Mathematics and Statistics, University of Melbourne, Parkville 3010, Australia. E-mail: s.mckinlay@pgrad.unimelb.edu.au.
}}

\date{}

\maketitle

\begin{abstract}
\noindent
The uniform law for sojourn times of processes with cyclically exchangeable increments is extended to the case of random fields, with general parameter sets, that possess a suitable invariance property.

\medskip


\end{abstract}

% ---------------------------------------------------------------------------------------------------------------------------------

\section{Introduction}

Among the most interesting and important problems on the pathwise behaviour of random processes is the one on the distribution of the time spent by the trajectory of the process under   a given (possibly random) line.  The simplest  question to be asked in that context is, of course, whether the trajectory of the process will cross that line at all. Historically, the first problem of that kind was perhaps the famous ballot theorem that effectively asserts that, for a simple symmetric random walk conditioned to be at point $-k<0$ at time~$n$, the probability that the trajectory will stay below zero on $(0,n]$ is $k/n$ (for a historical background and the development of ballot theorems, see~\cite{Tak}). The result is an immediate consequence of a simple invariance property of the random walk and can readily be extended to random walks with cyclicly exchangeable integer-valued jumps $\ge -1$ (see e.g.\ Theorem~7.6.1 in~\cite{Tak}).

It turned out that there are more interesting distributional results for sojourn times
that are also simple consequences of some invariance properties of the processes in
question. One of them is the uniform law for the sojourn time of the negative half-axis
in a simple symmetric random walk conditioned to be at zero at the terminal time, of
which the continuous time analog is the uniformity of the distribution of the time
spent below zero by the standard Brownian bridge on~$[0,1]$, the result going back
to~\cite{Levy1939} (that paper also contained the arc-sine law for the negative sojourn
time of the Brownian motion on~$[0,1]$). These laws were later generalised to skew
Brownian motion and skew Bessel processes, with  positive sojourn times that follow a
generalised arc-sine law~\cite{Lamperti,BaPiYo} (see also~\cite{Wa} for a
characterization result).

In the mid-1990s, the uniform laws for sojourn times were extended to suitable L\'evy bridges in~\cite{Fitzsimmons}, and necessary and sufficient conditions for such laws were provided for bridge processes with exchangeable increments in~\cite{Knight}. It appears that~\cite{Kallenberg1999}  was the first paper to note that the argument from~\cite{Knight}  also extends (with conditions) to any measurable bridge process with cyclically exchangeable increments (or,  equivalently, any measurable stationary periodic processes on $\reals_+$).

The nice invariance property exploited in~\cite{Knight}  allows  one to vastly simplify arguments proving some seemingly complex results (the reader may wish to compare the proofs in~\cite{Knight} with those in~\cite{Fitzsimmons}). This invariance property is discussed in \cite{Chaumont} and \cite{Yor2001}, and is due to the cyclic exchangeability of increments  which is defined as follows.

Let $X = \{X(t): t \in [0,1]\}$ be a measurable real-valued stochastic process. We say that $X$ has  cyclically exchangeable increments if, for any $u \in [0,1]$, the process
\[
X_u(t) := \left\{ \begin{array}{ll}
X(u+t)-X(u)+X(0) & \textrm{for $0 \le t <1-u$},\\
X(1)-X(u)+X(u+t-1)& \textrm{for $1-u \le t \le 1$},
\end{array}\right.
\]
has the same distribution as $X$. The simplest examples of such objects are L\'evy processes and bridges, and the uniform empirical processes.

Introduce the ``sojourn function"
\begin{equation*}
F(X,x) = \lambda \bigl(  t\in [0,1] :  X(t) \le x \bigr),
\quad x \in \reals,
\end{equation*}
where $\lambda$ is the Lebesgue measure on $\reals$. The value $F(X,x)$ is the (random)
duration of time the process $X$ spent below or at the level $x$. Then the following
result on the uniformity of the negative sojourn time holds true (see e.g.\
\cite{Kallenberg1999} or \cite{Yor2001}).

\begin{theorem}
  \label{thm1}
If $X$ has  cyclically exchangeable increments, $X(0)=X(1)=0$ and $F(X,\cdot)$ is continuous a.s., then $F(X,0)\sim U(0,1).$
\end{theorem}

In what follows, we will be interested in situations where $X(0)=X(1)$, in which case the definition of $X_u$ simplifies to
\[
X_u(t) :=
X(u+t \mbox{ (mod 1)})-X(u)+X(0)  , \quad t\in [0,1].
\]
In that case, it is more convenient  to view the process $X$ as given on the unit circle $\mathbb{S}^1$ using, say, the natural complex number parametrisation
\[
\widetilde{X}(e^{2\pi it}):= X (t), \quad t\in [0,1).
\]
Moreover, assuming integrability of $X$, observe that the cyclic exchangeability of increments now translates into the invariance, with respect to rotations of the parametric set $\mathbb{S}^1$, of the distribution of the process $\{\widetilde{X}^0 (u): u\in  \mathbb{S}^1\}$ given by
\[
\widetilde{X}^0 (e^{2\pi it}):= X (t) - \int_0^1 X(s) \, ds, \quad t\in [0,1).
\]
Now the point $t=0$ ceases to be special since, for any fixed $a\in \mathbb{S}^1,$ the distribution of the ``time" spent by $\widetilde{X}^0$ below the level $\widetilde{X}^0 (0)$ coincides with that of the time spent by $\widetilde{X}^0$ below the level $\widetilde{X}^0 (a)$, and in view of the result of Theorem~\ref{thm1}, that distribution is uniform.

The objective of the present note is to demonstrate how the above result can, in a natural (and rather elementary) way,  be extended to the random fields  setting where the parametric set $T$ of the field $\{X(t): t\in T\}$ is endowed with a finite measure $\mu$. We show that, for a  fixed $a\in T,$ provided that a certain invariance property is satisfied,  the $\mu$-measure of the subset of $T$ on which the values of $X$ do not exceed that of $X(a)$, is also uniformly distributed. In a sense, this result is a ``randomised" version of the well-known fact that, for a random variable $\xi$ with continuous distribution function $H$, one has $H (\xi)\sim U (0,1).$

\section{The main result}

Let $\pspace$ be a measure space with a finite measure and $\{X(t): t \in T\}$    a
real-valued measurable random field on a probability space $(\Omega, \mathscr{F}, \pr)$
with parameter set~$T$. Without losing generality, we will assume that $\mu(T) =1$.
Further, let
  \[
F_\mu (X, x) := \mu \bigl(  t\in T :  X(t) \le x \bigr),
  \quad x \in \reals,
  \]
%  for $a\in T,$
be the $\mu$-measure of the parameter   values $t$ for which $X(t)$ was at or below the
level~$x$.

Now suppose that $G$ is a family of measurable transformations of $T$, endowed with a
$\sigma$-algebra $\mathscr{G}$ containing cylinders with bases from $\mathscr{T},$
i.e.\ such that
\[
\{g\in G: g(t)\in B\}\in \mathscr{G}
 \quad\mbox{for any}\quad  t\in T,\ B\in \mathscr{T}.
\]
Denote by $X_g (\cdot) := X(g(\cdot))$ the composition of the random field $X$ and argument transformation $g\in G,$ and by $\mu_g (\cdot) := \mu (g(\cdot))$ a similar transformation of the measure~$\mu$.

A natural extension of Theorem~\ref{thm1} to this setting is given by

\begin{theorem}
 \label{thm2}
Suppose that $X_g\stackrel{d}{=}X$ and $\mu_g =\mu$ for any $g\in G,$  and that there
exists a probability measure  $\nu$ on $(G,\mathscr{G} )$ such that, for some fixed
$a\in T,$
\begin{equation}
\nu \bigl(  g\in G : g(a) \in B \bigr) =\mu (B), \quad B\in \mathscr{T}.
 \label{nu}
\end{equation}
If $F_\mu (X,\cdot)$ is continuous a.s.\ then     $F_\mu (X,X(a))\sim U(0,1).$
\end{theorem}

\begin{remark}[Remark 1.]
In the context of the situation discussed in the Introduction,  we had
$T=\mathbb{S}^1$, $\mu$ was the uniform distribution on~$\mathbb{S}^1$, the family $G$
consisted of all rotations of $\mathbb{S}^1$ (with $X_u$ corresponding to ``rotating"
the field $X$ by the  angle $2\pi u$), $\nu$ was the Haar measure on $G$. Clearly,
property~\eqref{nu} held  true for any $a\in T$.
\end{remark}

There is a discrete parameter set analog of the above theorem. Prior to stating it,
observe that, in the case of a finite parameter set $T,$ we can always take
$\mathscr{T}=2^T$ and that, assuming that $\mu$ is the uniform distribution on $T$, we
automatically have  $\mu_g =\mu$ for any transformation $g$ of $T$. Moreover, in that
case $G$ is also necessarily finite, so one can always take $\mathscr{G} = 2^G$.

The assertion has the following form.

\begin{theorem}
 \label{thm3}
Suppose that $X_g\stackrel{d}{=}X$ for any $g\in G,$ $\mbox{\rm card} (T)=N<\infty$,
$\mu (\{t\})=1/N$, $t\in T,$  and    that there exists a probability measure $\nu$ on
$G$ such that, for some fixed $a\in T,$ relation \eqref{nu} holds. If\/ $\pr (X(s) =
X(t))=0$  for any $s\neq t, $  then    $F_\mu (X,X(a))$ follows the uniform
distribution on $\{1,2,\ldots, N\}.$
\end{theorem}

The demonstration of the above result uses the same type of argument as that of
Theorem~\ref{thm2} and is actually even simpler, and so we will only present here the

\medskip

{\em Proof of Theorem~\ref{thm2}.} For $p\in (0,1),$ introduce
\begin{align*}
q_X (p)  &: = \sup\{x \in \reals : F_\mu (X, x) \le p\}, \\
Q_X (p) & : =\{t\in T: X(t) \le q_X (p)\}.
\end{align*}
By the continuity assumption,
\begin{equation}
 \label{fp}
F_\mu (X, q_X (p))\equiv \mu (Q_X (p))= p,
\end{equation}
and it is not hard to see that, for a fixed $t\in T,$ one has
\begin{equation}
 \label{in}
\{F_\mu (X,X(t)) \le p\} = \{t\in Q_X (p)\}.
\end{equation}
Observe also that, due to the invariance of $\mu$ under $g$, we have
\[
F_\mu (X_g , x) = \mu (t: X_g (t) \le x)
 = \mu_g (t: X_g (t) \le x)
 =  \mu  (s: X  (s) \le x) =F_\mu (X , x)
\]
for any $g\in G$, $x\in\reals$, and so $q_X (p) = q_{X_g} (p)$, $ p\in (0,1).$
Therefore, for any  $p\in (0,1),$ $g\in G$,
\begin{equation}
 \label{qqq}
  a \in Q_{X_g} (p) \quad\mbox{iff}\quad g(a) \in Q_{X } (p).
\end{equation}

Now consider the product probability space $(\Omega \times G, \mathscr{F} \otimes
\mathscr{G} , \pr \times \nu)$, and let  $\gamma$ be the second coordinate random
variable on the space, so that $\gamma$ is a random element of $G$ independent of $X$
and following the distribution~$\nu$. Then we clearly have $X_\gamma \stackrel{d}{=}
X,$ and so, for any $p\in (0,1)$,
\begin{align}
\pr\bigl(F_\mu (X,X(a)) \le p\bigr)
  &= \Pr\bigl(F_\mu (X_\gamma ,X_\gamma(a)) \le p\bigr)
   \notag \\
  &= \Pr \bigl(a \in Q_{X_\gamma} (p)\bigr)
  \label{pop} \\
    &= \Pr \bigl(\gamma(a) \in Q_{X } (p)\bigr)
  \label{quq}\\
&= \ex \Pr \bigl(\gamma(a) \in Q_X(p) \,|\,   X\bigr)
  \notag\\
&= \ex \,  \mu(Q_X(p))
 \label{pup} \\
&= p,
 \label{pap}
\end{align}
where \eqref{pop}, \eqref{quq}, \eqref{pup} and \eqref{pap} follow from \eqref{in},
\eqref{qqq}, \eqref{nu}  and~\eqref{fp}, respectively. The theorem is proved.

\section{Examples}

\begin{example}[Example 1.]
Let $U_1, \ldots, U_n$ be an i.i.d.\ sample of points uniformly distributed over
$\mathbb{S}^{d-1}$, $d\ge 2$, $K: \reals^{d} \rightarrow \reals$ be a  measurable
function, and set
\begin{equation}
 \label{XtEmp}
X(t) := \sum_{k=1}^n K(t-U_k), \quad t\in \mathbb{S}^{d-1}.
\end{equation}
Take $G:=SO(d)$, the special orthogonal group on $\reals^d$, and $\nu$ the Haar measure
on it. Then it is clear that, for a nice enough $K$, the field $X$ will satisfy the
conditions of Theorem~\ref{thm2} with $\mu$ being the uniform distribution on
$\mathbb{S}^{d-1}$ (due to symmetry).

To graphically illustrate this example, assume that $d=3$ and let
\begin{equation}
 \label{Ker}
K(x):=  1-(\|x\|/2)^{1/10} , \quad x\in \reals^{3},
\end{equation}
where $\|\cdot\|$ is the Euclidean norm. It is clear that the conditions of  Theorem~\ref{thm2} are met with the above choice of $K$.

One can think of $X$ as describing the relief (the surface  altitude relative to a
fixed reference level, i.e.\ the surface of a sphere) of a ``random planet" with $n$
mountains that have their summits at the sample points $U_k$ and possess symmetric
shape features (that are actually overlapping). A simulation of such a surface, with
$n=20$, is depicted on the left pane in Fig.~1.

\begin{figure}[h!]
 \label{fig}
    \centering
 \subfigure{\includegraphics[scale=0.53]{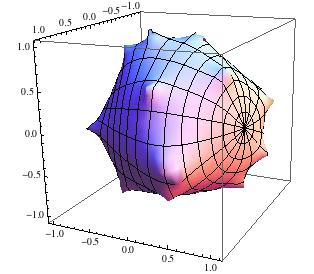}}
  \hspace{0.7cm}
  \subfigure{\includegraphics[scale=0.65]{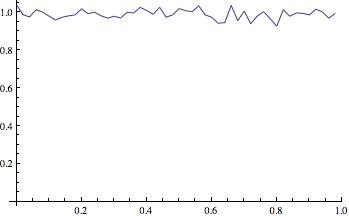}}
\caption{
 {\small
The left pane displays a sample realisation of random field \eqref{XtEmp} with kernel
\eqref{Ker}, $d=3$ and $n=20$. On the right pane, a scaled histogram of a simulated
sample of $10^5$ independent copies of~$F_\mu (X, X(a))$ for that field is shown, $a$
being the ``North pole".}
 }
\end{figure}

The right pane in Fig.~1 depicts the scaled histogram (with 50 bins) for an i.i.d.\ sample of $10^5$ copies of~$F_\mu (X, X(a))$ for that field, $a$ being the ``North pole" (corresponding to the centre of the ``spider web" seen on the left pane) on the ``coordinate sphere" $\mathbb{S}^{2}$. The values of $F_\mu (X, X(a))$ were evaluated using a simple antithetic Monte Carlo procedure: for each realisation of~$X$, $F_\mu (X, X(a))$ was estimated by the proportion of the $10^2$ points (half of which were chosen independently at random on $\mathbb{S}^{2}$, the other half being the opposite to these random points) at which the value of $X$ wouldn't exceed $X(a)$. The simulation confirms the assertion of Theorem~\ref{thm2} that the area on the coordinate sphere where the terrain is below the height of the relief at the North pole is
uniformly distributed.
\end{example}

\begin{example}[Example 2.]
Consider a rectangular array of random variables $\{X(t): t\in T\}$ with $T=
\{1,2,\ldots, m\}\times \{1,2,\ldots, n\}$ for some $n,m\in \naturals$, such that $\pr
(X(s) = X(t))=0$  for any $s\neq t.$ In other words, $X$ is a random $m\times n$ matrix
whose entries are distinct a.s. We will also assume that the distribution of the matrix
is invariant w.r.t.\ cyclic permutations of its rows and columns.

A real-life example of such a situation is provided by an $n\times n$ Latin square
experimental design, which is obtained by choosing at random a square from a given
isotopy class of $n\times n$ Latin squares (or, to simplify the procedure, starting
with a given $n\times n$ Latin square and then  cyclicly permuting its rows and columns
at random), and then assigning the treatment factor levels to the respectively labelled
elements of the array.

Introduce $G$ as the collection of all compositions of the form $g_r \circ g_c,$ where
$g_r$ is a cyclic permutation of the rows, and $g_c$ a cyclic permutation of the
columns of~$T$. Then clearly $X_g\stackrel{d}{=}X$ for any $g\in G$, and the uniform
distribution $\nu$ on $G$ will have the desired property \eqref{nu} for any fixed $a\in
T$ (recall that, in the discrete case, we assume that $\mu$ is uniform on $T$).

Now we see that all the conditions of Theorem~\ref{thm3} are met, and therefore, say,
the number of the entries of the matrix $X$ that don't exceed the first entry in its
first row (i.e.\ the random variable $ \mbox{\rm card} \bigl(  t\in T : X(t) \le
X((1,1)) \bigr)$ has the uniform distribution on $\{1,2, \ldots, nm\}.$
\end{example}

\end{document}